\documentclass[12pt,nodate]{article}
\usepackage{amsmath,amssymb,amsthm}
\date{}
\begin{document}

\newtheorem{theorem}{Theorem}
\newtheorem{lemma}{Lemma}
\newtheorem{proposition}{Proposition}
\newtheorem{definition}{Definition}
\newtheorem{corollary}{Corollary}
\newtheorem{notation}{Notation}

\newcommand{\cF}{\ensuremath{\mathcal{F}}}
\newcommand{\PRF}{\ensuremath{\mathbb{R}(PR(\mathcal{F}))}}
\newcommand{\C}{\ensuremath{\mathbb{C}}}
\newcommand{\Cn}{\ensuremath{\mathbb{C}^{n}}}
\newcommand{\PF}{\ensuremath{PR(\mathcal{F})}}
\newcommand{\wdF}{\ensuremath{\mathbf{w} \in dom(F)}}
\newcommand{\dFi}{\ensuremath{\frac{\partial F}{\partial z_{i}}}}
\newcommand{\bz}{\ensuremath{\mathbf{z}}}
\newcommand{\R}{\ensuremath{\mathbb{R}}}
\newcommand{\dFsn}{\ensuremath{\frac{\partial F}{\partial z_{n+1}}}}
\newcommand{\bw}{\ensuremath{\mathbf{w}}}
\newcommand{\ba}{\ensuremath{\mathbf{a}}}
\newcommand{\cFt}{\ensuremath{\tilde{\mathcal{F}}}}
\newcommand{\awd}{\ensuremath{\langle w_{1} , \ldots , w_{n} \rangle}}
\newcommand{\idFsn}{\ensuremath{\frac{\partial^{i}F}{\partial z_{n+1}^{i}}}} 
\newcommand{\pidFsn}{\ensuremath{\frac{\partial^{i-1}F}
{\partial z_{n+1}^{i-1}}}}
\newcommand{\dHssn}{\ensuremath{\frac{\partial H}{\partial z_{n+2}}}}
\newcommand{\tF}{\ensuremath{\tilde{\mathcal{F}}}}
\newcommand{\dfi}{\ensuremath{\frac{\partial f}{\partial z_{i}}}}
\newcommand{\dGij}{\ensuremath{\frac{\partial G_{i}}{\partial z_{j}}}}
\newcommand{\tD}{\ensuremath{\tilde{D}}}
\newcommand{\lm}{\ensuremath{[\lambda : \mu ]}}
\newcommand{\bx}{\ensuremath{\mathbf{x}}}
\newcommand{\by}{\ensuremath{\mathbf{y}}}
\newcommand{\mx}{\ensuremath{\mu(x)}}
\newcommand{\lx}{\ensuremath{\lambda(x)}}
\newcommand{\ly}{\ensuremath{\lambda(y)}}
\newcommand{\si}{\ensuremath{\sqrt{-1}}}
\newcommand{\my}{\ensuremath{\mu(y)}}
\newcommand{\bxy}{\ensuremath{(\mathbf{x}, \mathbf{y})}}
\newcommand{\abxy}{\ensuremath{\langle \mathbf{x}, \mathbf{y} \rangle}}
\newcommand{\pbxy}{\ensuremath{\mathbf{x} + \sqrt{-1} \mathbf{y}}}
\newcommand{\duxi}{\ensuremath{\frac{\partial u}{\partial x_{i}}}}
\newcommand{\duyi}{\ensuremath{\frac{\partial u}{\partial y_{i}}}}
\newcommand{\mxi}{\ensuremath{\mu(x_i )}}
\newcommand{\myi}{\ensuremath{\mu(y_i )}}
\newcommand{\lxi}{\ensuremath{\lambda(x_i)}}
\newcommand{\lyi}{\ensuremath{\lambda(y_i)}}
\newcommand{\cE}{\ensuremath{\mathcal{E}}}
\newcommand{\Rn}{\ensuremath{\mathbb{R}^{n}}}

\title{Analytic continuation and Zilber's quasiminimality conjecture}
      
\author{A.J. Wilkie\\
    Mathematical Institute\\
    University of Oxford\\}

\maketitle
\date{}

\vspace{5mm}

\hspace{53mm} $\mathbf{Abstract}$

\vspace{3mm}

In this article, which is dedicated to my friend and colleague Boris Zilber on the occasion of his 75th birthday, I put forward a strategy for proving his quasiminimality conjecture for the complex exponential field. That is, for showing that every subset of $\C$ definable in the expansion of the complex field by the complex exponential function is either countable or cocountable. In fact the strategy applies to any expansion of the complex field by a countable set of entire functions (in any number of variables) and is based on a certain property-an analytic continuation property-of the o-minimal structure obtained by expanding the ordered field of real numbers by the restrictions to compact boxes of the real and imaginary parts of the functions in the given set.

In a final section I discuss briefly the (rather limited) extent of our unconditional knowledge in the area.

\vspace{3mm}

\hspace{48mm} ........................

\vspace{3mm}

After some reflection we, Boris and I, agreed that it was in July 1993 that he first asked me whether I had thought about the model theory of the complex exponential field. The occasion was Logic Colloquium '93, the European Summer Meeting of the Association for Symbolic Logic, which was held that year at Keele University in the UK. We were both invited to give plenary lectures. Boris spoke about his results on model-theoretic dimensions in complex (and ultrametric) analysis and I gave a talk on recent work with Angus (eventually published in the paper [MW]) concerning the decidability of the real exponential field, which was still my main concern. I certainly hadn't even considered the complex exponential. In fact the idea seemed absurd; the set of integers was definable  and so as far as I was concerned it was not a tame structure. But Boris, then as now, had considerably more imagination, and deeper insights than I did into potentialy good model-theoretic behaviour of familiar mathematical structures. If the set of reals was also definable then, he agreed, the situation would indeed be hopeless. But what if it wasn't? In fact, could it not be the case that every definable subset of \C\ was either countable or co-countable? That remark had a profound effect on my mathematical life and was the motivational force for much of my research since that moment almost thirty years ago. It still is.

\vspace{3mm}

I soon realised that I could make progress on the problem-the \emph{quasiminimality} problem as Boris called it-only through o-minimal theory and hard analysis: my knowledge of abstract stability theory and its generalizations was simply not strong enough, and certainly not up to a level to be able to follow and make a meaningful contribution to Boris's ingenious, beautiful and eventually highly successful construction  of a ``pseudo-exponential'' field using a combination of Hrushovski's predimensions and Shelah's theory of excellent classes. So my own approach started with the observation that the complex exponential function is definable in an o-minimal structure (when considered as a function from $\mathbb{R}^{2}$ to $\mathbb{R}^{2}$) provided that the imaginary part of its argument is restricted to a bounded interval. (One may take this structure to be, for example, the real ordered field expanded by both the real exponential function and by the sine function restricted to the interval $[0, \pi]$.) After several years I came up with about forty pages of notes purporting to contain, amongst other things, a proof that the complex exponential field was quasiminimal and I intended to present a sketch of the argument at an Oberwolfach meeting in July 2004. 

\vspace{3mm}

Most readers of this article will know that Boris' conjecture is still unresolved and hence that there must have been a mistake, as indeed there was. But three things did come out of those notes. Firstly, I had established (correctly) a positive solution of Schanuel's conjecture for ``generic'' finite sequences of real numbers, and this rescued my Oberwolfach talk. ( The result was eventually generalized and published in the paper [BKW] written jointly with two of the editors of this volume.) Secondly,  in order to cope with the fact that the $2\pi i \mathbb{Z}$-periodicity of the complex exponential function is at complete variance with the whole ethos of o-minimality, I was forced to investigate whether, and if so how, points with  integer coordinates could lie in sets defined by algebro-exponential equations, and thereby obtain a measure of ``non-tameness''. I certainly did not want to assume the full version of Schanuel's Conjecture (which would have settled this particular issue), but if it could be shown that either there were few such points or, if not, then at least one such point behaved generically then the ``generic Schanuel theorem'' just mentioned could be used and this might be sufficient to control the periodicity. It didn't work but, and I apologize for a small personal digression at this point,  the idea of quantitatively investigating the occurrence of points with integer coordinates in sets defined as above did appeal to me. In fact, why not do the same in the framework of general o-minimal structures? So I temporarily abandoned work on the quasiminimalty problem and initiated the project of counting integer points lying in sets definable in a fixed, but arbitrary, o-minimal structure. I had limited success at the time, only managing to prove a result in the case of one dimensional sets. The paper ([W1]) was published in 2004 and came to the attention of Jonathan Pila who had been working on  similar issues. He had obtained bounds of the same general nature as mine for rational points lying in low dimensional, globally subanalytic sets (i.e. sets definable in the o-minimal structure $\R _{an}$), but was finding it difficult to put his arguments into a general setting that would smoothly facilitate an induction on dimension. Jonathan and I met soon after and agreed that the setting had to be o-minimality. Our point counting theorem was published in 2006 and this is how, at least from my point of view, that result came about: my motivation for studying integer points in o-minimally definable sets was, apart from the fun of it, completely motivated by Boris' quasiminimality problem.  
I had no inkling of how, over the following fifteen years, the result would be applied, with huge success, by Jonathan and many others to problems in diophantine geometry. 

Serendipity? Not exactly. For while Boris, as far as I know, did not foresee such applications, his conjectures, often fearless but always with sound intellectual bases, have been an inspirational source of  research throughout our community even when, and perhaps especially when, that research takes unexpected directions.

\vspace{3mm}

I mentioned that three things came out of my 2004 notes. The third occupies the remainder of this paper. It concerns my strategy, which I still believe to be plausible, for proving the quasiminimality conjecture for the complex exponential field and, possibly, for many other expansions of the complex field by entire functions. I am very grateful to the editors for giving me this opportunity to explain the strategy despite the fact that, at the time of writing, it has not resulted in any definite results. The main theorems do establish the quasiminimality of a certain class of structures expanding the complex field but is conditional on a property of locally definable holomorphic functions, namely that they have analytic continuations along ``generic'' paths (i.e. those avoiding \emph{obvious} singularities). The first theorem is precisely this, while the second provides a  criterion for the continuation in purely complex analytic terms that avoid notions of general definability.

\vspace{7mm}

So let $\widetilde{\R}$ be a fixed o-minimal expansion of the real ordered field \\ $\overline{\R} := \langle \R ; < , + , \cdot \rangle$. Following [PS1, PS2] we say that a complex valued function $F$ of the $n$ complex variables $z_1, \ldots, z_n$ is \emph{definable} if the real and imaginary parts of $F$ are definable in 
$\widetilde{\R}$ (\emph{without} parameters) when considered as functions of the $2n$ real variables $x_1,y_1,\ldots,x_n,y_n$, where $z_j=x_j+ i y_j$ for $j=1.\ldots,n$. If $F$ is holomorphic, its domain will always be assumed to be, without further mention, a connected, open subset (a fortiori a \emph{definable} subset if $F$ is definable) of \Cn . We write $dom(F)$ for the domain of $F$.

\vspace{3mm}

Peterzil and Starchenko develop complex analysis in this definable context, but in fact many of the subtleties of their work will not be needed here since we will only be working with the standard structure $\widetilde{\R}$. However, one of their results is still worth mentioning, namely that if $\{F_t:t\in \R^d \}$ is a definable family of $n$-variable, holomorphic functions, then there exists $N$ such that for all $t \in \R^d$ and all $a \in dom(F_t )$, if $\frac{\partial^{\alpha}F}{\partial z^{\alpha}}$ vanishes at $a$ for each $\alpha \in \mathbb{N}^n$ with $|\alpha | \leq N$, then $F_t$ is identically zero. (This is, of course, obvious if $\widetilde{\R}$ is polynomially bounded, but not so clear otherwise. For example, and of particular relevance to Zilber's problem, the $F_t$'s could range over polynomials in
$z_1, \ldots, z_n, \log z_1,\ldots,\log z_n$ of some fixed degree (and with suitable simply connected domains).)

\vspace{3mm}

Anyway, returning to definitions we consider the closure operator $LD(\cdot)$ - ``locally definable from'' - where we specify that for all $n \geq 1$ and for all  $a_1, \ldots, a_{n+1} \in \C$,

\vspace{2mm}

\noindent 
(1) \hspace{2mm} $a_{n+1} \in LD(a_1, \ldots, a_n )$ \hspace{1mm} :$\Longleftrightarrow$
\hspace{1mm} $F(a_1, \ldots, a_n ) 
= a_{n+1}$ \hspace{2mm}for some definable, holomorphic function $F$ with $\langle a_1, \ldots, a_n \rangle \in dom(F)$.

\vspace{2mm}

It is straightforward to show that $LD$ is a pregeometry and, further, that for all $b_1, \ldots, b_m \in \C$,

\vspace{2mm}

\noindent
(2)\hspace{2mm} $b_1, \ldots, b_m$ are $LD$-independent if and only if whenever $G$ is a nonzero, definable, holomorphic function with $\langle b_1, \ldots, b_m \rangle \in dom(G)$, then $G(b_1, \ldots, b_m) \neq 0$.

\vspace{2mm}

\noindent
(See [W3] for details.)

\vspace{3mm}

Strictly speaking $LD(X)$ should be specified for all subsets $X$ of \C\  and this is done in the usual way by taking $LD(X)$ to be the union of all $LD(X_0)$'s as $X_0$ ranges over finite subsets of $X$. There is a small issue concerning the $LD$-closure of the empty set which is taken to be $\{ s+ i t : s,t \in \R,$ both definable in  $\widetilde{\mathbb{R}} \}$. One easily checks (using the Cauchy-Riemann equations) that if $G$ is a nonzero, definable, holomorphic
function with $s+ i t \in dom(G)$ and  $G(s+ i t) = 0$, then both $s$ and $t$ are indeed definable in $\widetilde{\mathbb{R}}$, so this choice of $LD(\emptyset )$ is consistent with (2) above. However, the relationship between the $LD$-dimension of an arbitrary finite subset of \C\ and the dimension of the corresponding set of real and imaginary parts (for the usual pregeometry of the o-minimal structure $\widetilde{\mathbb{R}}$) is more complicated and is resolved in section 4 of [W3].

\vspace{3mm}

\noindent
$\mathbf{Examples}$

\vspace{2mm}

\noindent 
(a) If we take $\widetilde{\mathbb{R}}$ to be just $\overline{\mathbb{R}}$, then $LD$ is just algebraic closure (over $\mathbb{Q}$) in the field \C\ . More generally, if we expand $\overline{\mathbb{R}}$ by a set $A$ of constants, then $LD$ is algebraic closure over the subfield $\mathbb{Q}(A)$ of \C\ . (Two remarks: firstly, I do not find this at all obvious (the proof can be found in [W2]), and secondly, recall that our holomorphic functions are allowed to be definable in $\widetilde{\mathbb{R}}$ so, for example we can distinguish (uniformly in parameters) between different roots of a polynomial.)

\vspace{2mm}

\noindent 
(b) Let $\widetilde{\mathbb{R}} = \R_{exp} := \langle \overline{\mathbb{R}}, exp \rangle$ where $exp: \R \rightarrow \R : x \mapsto e^x$. Then, by a result of Bianconi (see [B]), $LD$ is still just algebraic closure (but this time it's over the minimal model $K \preccurlyeq \mathbb{R}_{exp}$ of the theory of
$\mathbb{R}_{exp}$). Of course, things are different, and of great relevance to Zilber's problem, if we expand $\R_{exp}$ by the restricted sine function and we will be discussing this situation later.

\vspace{2mm}

\noindent
(c) For $\widetilde{\mathbb{R}} = \R_{an}$ (where the definable sets are the globally subanalytic sets, see for example, [DD]), $LD$ is trivial: $LD(X) = LD(\emptyset ) = \C$ \hspace{1mm} for all $X \subseteq \C$ simply because all $r \in \R$ are definable.

\vspace{2mm}

In view of example (c) above we now assume that the language of $\widetilde{\mathbb{R}}$ is countable so, in particular, $LD(\emptyset )$ is countable. We shall require rather more:

\vspace{3mm}

\noindent
$\mathbf{Lemma}$ (Existence of generic lines)

\vspace{2mm}

Let $A$ be a countable subset of \R . Then there exists $\omega \in \C$ (in fact an uncountable, dense set of them) such that for all continuous functions $\phi : [0, 1] \mapsto \C$ definable in the structure $\langle \widetilde{\mathbb{R}}, s \rangle_{s \in A}$ and satisfying $\phi (0) \neq 0$, we have that $\phi (t) \neq t\omega$ for all $t \in [0, 1]$.

\vspace{1mm}

\noindent
$\mathbf{Proof}$

\vspace{1mm}

For each such $\phi$ let 

\vspace{2mm}

$S_{\phi} := \{ \omega \in \C$ : for all  $t \in [0, 1], \phi (t) \neq t\omega \}$.

\vspace{2mm}

It clearly follows from the continuity of $\phi$ and the compactness of the closed unit interval  that $S_{\phi}$ is an open subset of \C . So by the Baire Category Theorem we will be done if we show that $S_{\phi}$ is a dense subset of \C\  (because there are only countably many $\phi$'s).

So let $\omega_{0}$ be an arbitrary complex number and let $\epsilon > 0$. Let $\Delta$ be the open disc in \C\ centered at $\omega_0$ and of radius $\epsilon$ and suppose, for a contradiction, that $\Delta \subseteq \C \smallsetminus S_{\phi}$, i.e. that for all $\omega \in \Delta$, there exists $t \in [0, 1]$ such that $\phi (t) = t\omega $. Then, by definable choice and the non-decreasing of dimension under definable, injective maps (working in the o-minimal structure $\langle \widetilde{\mathbb{R}}, s \rangle_{s \in A}$), there exist $\omega_1 , \omega_2 \in \Delta$ with $\omega_1 \neq \omega_2$ such that for some $t_0 \in [0, 1]$, $\phi (t_0 ) = t_0 \omega_1 $ and $\phi (t_0 ) = t_0 \omega_2 $. This is absurd unless $t_0 = 0$; but this is also ruled out since $\phi (0) \neq 0$. $\Box$

\vspace{3mm}

We will be considering analytic continuations of definable functions along generic paths in \Cn . In fact, we will only need to consider linear paths: for $a, \omega \in \Cn$ ($n \geq 1$), define the map $\lambda_{a, \omega} : \C \mapsto \Cn$ by $\lambda_{a, \omega} (z) := a + z\omega $ (for $z \in \C$). We say that  $\lambda_{a, \omega}$ is \emph{generic} on a set $T \subseteq \C$ if $\lambda_{a, \omega} (t)$ is a generic $n$-tuple for each $t \in T$, i.e. if $a_1 + t\omega_1  , \ldots , a_n + t\omega_n $ are $LD$-independent complex numbers for each $t \in T$ (where $a = \langle a_1 , \ldots , a_n \rangle$ and 
$\omega = \langle \omega_1 , \ldots , \omega_n \rangle$). The set $T$ will, almost always, be the interval $[0, 1]$ so that for any $a, \omega$ we have that $\lambda_{a, \omega}(0) = a$.

\vspace{1mm}

We now come to our main definitions.

\newpage

\noindent
$\mathbf{Definition}$ $\mathbf{1}$

\vspace{1mm}

We say that the structure $\widetilde{\R}$ has the \emph{Analytic Continuation Property} (ACP) if for all $LD$-independent $a_1 , \ldots , a_n \in \C$, all definable, holomorphic functions $F$ with $a= \langle a_1 , \ldots , a_n \rangle \in dom(F)$ and all $\omega \in \Cn$ with $\lambda_{a, \omega}$ generic on $[0, 1]$, there exists a definable, holomorphic function $G$ with  $\lambda_{a, \omega}([0, 1]) \subseteq dom(G) \subseteq \Cn$ such that $G(a)  = F(a)$. (And hence, by (2), $G(\lambda_{a, \omega} (z)) = F(\lambda_{a, \omega} (z))$ for all $z \in \C$ such that $\lambda_{a, \omega} (z)$ lies in the connected component of $dom(F) \cap dom(G)$ containing the point $a$. So the function $G\circ \lambda_{a, \omega}$  analytically continues, in the usual sense,  the function $F\circ \lambda_{a, \omega}$ (restricted to a sufficiently small open neighbourhood of $0 \in \C$) to an open set containing the interval $[0, 1]$.)

\vspace{3mm}

\noindent 
$\mathbf{Definition}$ $\mathbf{2}$

\vspace{1mm}

Let $1 \leq l \leq n$ and let $M \subseteq \Cn$. Then we say that $M$ is an $l$-\emph{dimensional, locally definable, complex submanifold of} \Cn\ (or just an $l$-\emph{manifold} for short) if

\vspace{1mm}

\noindent
(a) $M$ is a closed subset of \Cn  , and 

\vspace{1mm}

\noindent
(b) for all $a \in M$, there exist a definable open set $W$ with $a \in W \subseteq \Cn$ and a holomorphic, definable map $G = \langle G_1 , \ldots , G_{n-l} \rangle : W \mapsto \mathbb{C}^{n-l}$ such that $a$ is a non-singular point of the zero set of $G$ (i.e. $G(a) = 0$ and the vectors $\langle \frac{\partial G_j}{\partial z_1} (a) , \ldots , \frac{\partial G_j}{\partial z_n} (a) \rangle$ (for $1 \leq j \leq n-l$) are linearly independent over \C\ ) and, further, $M \cap W = Z_{reg} (G)$ where $Z_{reg} (G)$ denotes the set of non-singular points of the zero set of $G$.

\vspace{3mm}

\noindent
$\mathbf{Example}$

\vspace{1mm}

If $\widetilde{\R} = \langle \overline{\R} , exp$$\upharpoonright[0, 1]$, $sin$$\upharpoonright$$[0, 2\pi ] \rangle$, then the graph of the complex exponential function $\Gamma_{exp} := \{ \langle z, e^z \rangle : z \in \C \}$ is a $1$-dimensional, locally definable, complex submanifold of $\C^2$.

\vspace{3mm}

\noindent
$\mathbf{Definition}$ $\mathbf{3}$

\vspace{1mm}

The structure $\widetilde{\C}$ is defined to be the expansion of the complex field by all $l$-dimensional, locally definable, complex submanifolds of \Cn\  (for all $l$, $n$ with $1 \leq l < n$).

\vspace{2mm}

\noindent
$\mathbf{Remark}$

\vspace{1mm}

With $\widetilde{\R}$ as in the example above we see that $\widetilde{\C}$ is an expansion of $\C_{exp}$ (:= the complex field expanded by the complex exponential function). I do not know if it is a proper expansion. For example, it is clear that every connected component of an $l$-manifold is also an $l$-manifold, but it seems to me (and, in fact, to Zilber) to be perfectly possible that some $l$-manifold, $M$ say, is definable in $\C_{exp}$ but that some connected component of $M$ is not. 

\vspace{3mm}

\noindent$\mathbf{Theorem}$ $\mathbf{1}$

\vspace{1mm}

Suppose that $\widetilde{\R}$ has the ACP. Then the structure $\widetilde{\C}$ is quasi-minimal. That is, for any subset $S$ of $\C$ which is definable in the language $\mathcal{L} (\widetilde{\C})$ of $\widetilde{\C}$, and we do allow parameters here, we have that either $S$ is countable or its complement $\C \setminus S$ is countable. In fact, the same is true for sets $S$ defined by a formula of the infinitary language $\mathcal{L} (\widetilde{\C})_{\infty , \omega}$ provided, of course, that the formula contains only countably many parameters.

\vspace{1mm}

\noindent
$\mathbf{Proof}$

\vspace{1mm}

We first show that if $u$ and $v$ are elements of $\C \setminus LD(\emptyset )$ then there exists a back-and-forth system (for the structure $\widetilde{\C}$) containing the pair $\langle u, v \rangle$.

Then by a classical result of Karp ([K]) the quasiminimality condition for parameter-free formulas follows from this since it implies that if $S \subseteq \C$ is a parameter-free $\mathcal{L} (\tilde{\C})_{\infty , \omega }$-definable set, then either $S \subseteq LD(\emptyset)$ or else $\C \setminus LD(\emptyset ) \subseteq S$, and $LD(\emptyset)$ is countable.

So suppose that $u, v \in \C \setminus LD(\emptyset)$. We may assume that $\lambda_{u, v-u}$ is generic on $[0, 1]$.
For otherwise, by the countability of $LD(\emptyset )$, there exists some $w \in \C \setminus LD(\emptyset)$ such that both $\lambda_{u, w-u}$ and 
 $\lambda_{w, v-w}$ are generic on $[0, 1]$ and we prove the result for the pair $u, w$ and for the pair $w, v$.
 
 We now set up a back-and-forth argument. 
 
 \vspace{2mm}
 
 For $n \geq 1$ and $a = \langle a_1 , \ldots , a_n \rangle , b = \langle b_1 , \ldots , b_n \rangle \in \Cn$ we write $a \sim_n b$ if
 
 \vspace{2mm}
 
 \noindent
 $(i)_n$ $ a_1 = u$ and $b_1 = v$ ;

\vspace{2mm}

\noindent
$(ii)_n$ for some $m$ with $1 \leq m \leq n$ and some $i_1 , \ldots , i_m$ with $1 \leq i_1 < \cdots < i_m \leq n$ we have that $\lambda_{a' , b'-a'}$ is generic on $[0, 1]$, where $a' := \langle a_{i_1} , \ldots , a_{i_m} \rangle$
 and  $b' := \langle b_{i_1} , \ldots , b_{i_m} \rangle$;
 
\vspace{2mm}

\noindent
$(iii)_n$ there exists a definable connected, open set $V \subseteq \mathbb{C}^{m}$ with $\lambda_{a' , b'-a'}([0, 1]) \subseteq V$ and, for each $j = 1, \ldots , n$, a definable, holomorphic function $F_j : V \mapsto \C$ such that
\hspace{2mm} $F_j (\lambda_{a' , b'-a'}(0)) =F_j (a_{i_1} , \ldots , a_{i_m}) = a_j$ and 
\hspace{2mm} $F_j (\lambda_{a' , b'-a'}(1)) =F_j (b_{i_1} , \ldots , b_{i_m}) = b_j$.

\vspace{3mm}

We clearly have $\langle u \rangle \sim_1 \langle v \rangle$ (where, in this case, $n = m = 1$, $V = \C$ and $F_1$ is the identity function on \C\ ) .

\vspace{3mm}

In order to establish the back-and-forth property, suppose that $n \geq 1$ and that $a \sim_n b$ as above, with $m$, $a'$, $b'$ as in $(ii)_n$. We write $\lambda$ for $\lambda_{a', b'-a'}$.

\vspace{3mm}

So let $a_{n+1} \in \C$. We must find $b_{n+1} \in \C$ such that $\langle a, a_{n+1} \rangle \sim_{n+1} \langle b, b_{n+1} \rangle$. There are two cases.

\vspace{3mm}

\noindent
\emph{Case 1}. $a_{n+1} \notin LD(a_{i_1} , \ldots , a_{i_m})$. 

\vspace{2mm}

Let $A$ be a finite subset of \R\ containing the real and imaginary parts of
$a_{i_1} , \ldots , a_{i_m} , a_{n+1} , b_{i_1} , \ldots , b_{i_m}$. Apply the generic lines lemma to obtain some $\omega \in \C$ such that for all continuous $\phi : [0, 1] \rightarrow \C$ definable in the structure $\langle \tilde{\R}, c \rangle_{c \in A}$ with $\phi (0) \neq 0$, we have $\phi (t) \neq t\omega $ for all $t \in [0, 1]$ and, further,
such that $\omega$ does not lie in the (countable) set $LD(b_{i_1} , \ldots , b_{i_m} , a_{n+1})$.
Let $b_{n+1} := a_{n+1} + \omega$. We show that $\langle a, a_{n+1} \rangle \sim  \langle b, b_{n+1} \rangle$. 
Now $(i)_{n+1}$ is obvious and for $(ii)_{n+1}$, $(iii)_{n+1}$ we replace $m$ by $m+1$ ($\leq n+1$) and let $i_{m+1} := n+1$. We also replace $a'$ by $\langle a', a_{n+1} \rangle$ and $b'$ by $\langle b', b_{n+1} \rangle$ so that $\lambda$ ($=  \lambda_{a', b'-a'}$) becomes $\lambda^* : \C \mapsto \C^{m+1}$ given by $\lambda^* (z) := \langle \lambda (z), a_{n+1} + z\omega \rangle$. We must show that $\lambda^*$ is generic on $[0, 1]$. So suppose, for a contradiction, that for some $t_0 \in [0, 1]$ the $(m+1)$-tuple $\lambda^* (t_0 )$ has $LD$-dimension $< m+1$. Now since $\lambda (t_0 )$ has $LD$-dimension $m$ (because $\lambda$ is generic on $[0, 1]$) we must have $a_{n+1} + t_0 \omega \in LD(\lambda (t_0 ))$. So there exists a definable, holomorphic function $F$ with $\lambda (t_0 ) \in dom(F)$ and such that

\vspace{3mm}

\noindent
(*)   \hspace{20mm} $F(\lambda (t_0 )) = a_{n+1} + t_0 \omega $.

\vspace{3mm}

Now we cannot have $t_0 = 0$ for then we would have $a' = \lambda (0) \in dom(F)$ and  $F(a' ) = a_{n+1}$ which contradicts the hypothesis of Case 1. Also, $t_0 \neq 1$ for otherwise we would have that $\omega = F(\lambda (1)) - a_{n+1} =  F(b') - a_{n+1} \in LD( b_{i_1} , \ldots , b_{i_m} , a_{n+1})$.
So $0 < t_0 < 1$ and since $dom(F)$ is open and $\lambda$ is continuous, there
exist rationals $q_1$, $q_2$ with $0 < q_1 < t_0 < q_2 < 1$ such that $\lambda ([q_1 ,q_2 ]) \subseteq dom(F)$.

Define $\phi : [0, 1] \rightarrow \C$ to be the continuous function
taking the value $-a_{n+1} + F(\lambda (t))$ for $q_1 \leq t \leq q_2$  and which is linear on each of the intervals $[0, q_1 ]$, $[q_2 ,1]$ with, say, $\phi (0) = 1$ and $\phi (1) = 0$.

Then $\phi$ is definable in the structure $\langle \widetilde{\R}, c \rangle_{c \in A}$ because $\lambda$ is, $a_{n+1}$ is, and $F$ is definable in $\widetilde{\R}$ without parameters. Further, $\phi (0) \neq 0$ so by the construction of $\omega$ we have that $\phi (t) \neq t\omega $ for all $t \in [0, 1]$. In particular, $\phi (t_0 ) \neq t_0$, i.e. $-a_{n+1} + F(\lambda (t_0 )) \neq t_0 \omega$, which contradicts (*) and establishes $(ii)_{n+1}$.

\vspace{2mm}
 
 As for $(iii)_{n+1}$, we take our new $V$ to be $V \times \C$ and our new $F_j$'s, call them $F^{*}_j$ (for $j = 1, \ldots , n+1$), to be specified  (for $\langle z_1 , \ldots , z_{m+1} \rangle \in V \times \C$) by setting $F_j^{*} (z_1 , \ldots , z_{m+1}) := F_j (z_1 , \ldots , z_{m})$ if $1 \leq j \leq n$ and $F_{n+1}^{*} (z_1 , \ldots , z_{m+1}) := z_{m+1}$.

Then the required conditions for $(iii)_{n+1}$ carry over from $(iii)_n$ for 
$j = 1, \ldots , n$, and for $j=n+1$ we have, for each $t \in [0, 1]$,

\vspace{3mm}

\hspace{4mm} $F_{n+1}^{*} (\lambda _{\langle a' , a_{n+1} \rangle, \langle b' , b_{n+1}   \rangle - \langle a' , a_{n+1} \rangle} (t)) = a_{n+1} + t(b_{n+1} - a_{n+1})$

 \vspace{3mm}

 \noindent
which takes the value $a_{n+1}$ for $t=0$ and $b_{n+1}$ for $t=1$. So we have that $\langle a, a_{n+1} \rangle \sim_{n+1} \langle b, b_{n+1} \rangle$ as required.
 
 \vspace{3mm}

 \noindent
 \emph{Case 2}. $a_{n+1} \in LD(a_{i_1} , \ldots , a_{i_m})$. 

\vspace{2mm}
In this case there exists a definable, holomorphic function $F$ such that $a' \in dom(F)$ and $F(a') = a_{n+1}$. Apply the ACP to obtain a definable, holomorphic function $G$ with $\lambda ([0, 1]) \subseteq dom(G) \subseteq \C^m$ satisfying $G(a') = F(a')$, i.e. $G(a') = a_{n+1}$. Now with $V \subseteq \C^m$ as given by $(iii)_{n}$, note that  $\lambda ([0, 1]) \subseteq V \cap dom(G)$. Let $U$ be a definable, connected open subset of $V \cap dom{G}$ such that $\lambda ([0, 1]) \subseteq U$. (Clearly such a $U$ exists and may be taken, for example, to be a certain finite union of polydiscs with Gaussian rational centres and rational radii.) 

We now take the same $m$ as in $(ii)_{n}$ so that $(i)_{n+1}$ and $(ii)_{n+1}$ are obviously satisfied. For $(iii)_{n+1}$ we take the $F_j$'s as given by  $(iii)_{n}$ for $j = 1, \ldots , n$ and restrict them to the set $U$. For $F_{n+1}$ we take $G$ restricted to $U$ so that $F_{n+1} (\lambda (0)) =  F_{n+1} (a') = G(a') = a_{n+1}$. Finally, taking $b_{n+1} := F_{n+1} (\lambda (1))$ completes the construction in Case 2.
 
 \vspace{5mm}
 
 Of course we also need to consider the ``back'' case-where we take some $b_{n+1} \in \C$ and have to find some $a_{n+1} \in \C$ satisfying $\langle a, a_{n+1} \rangle \sim \langle b, b_{n+1} \rangle$. But this follows in exactly the same way upon noting that for any $c, d \in \C^l$ we have $\lambda_{c, d-c}$ is generic on $[0, 1]$ if and only if $\lambda_{d, c-d}$ is (because the ranges on $[0, 1]$ are the same).
 
 \vspace{2mm}
 
 So, our system $\{ \langle a, b \rangle : n \geq 1, a, b \in \Cn$ and $a \sim_n b \}$ has the back-and-forth property.
 
 We must now show that atomic formulas of $\mathcal{L}(\tilde{\C})$  are preserved from $a$ to $b$ whenever $a \sim_n b$.
 
 So let $1 \leq l \leq n$ and suppose that $M \subseteq \Cn$ is an $l$-manifold. Suppose $a \in M$ and $b \in \Cn$ are such that $a \sim_n b$. We must show that $b \in M$ (and similarly for $a$ and $b$ interchanged, for which the proof is the same). Write $a = \langle a_1 , \ldots , a_n \rangle$ and $b = \langle b_1 , \ldots , b_n \rangle$. Let $m$, $1 \leq i_1 < \cdots < i_m \leq n$, $V$ and the $F_j$'s be as in $(ii)_n$ and $(iii)_n$, and write $F$ for the map $\langle F_1 , \dots , F_n \rangle : V \rightarrow \Cn$.
 
 Define $T := \{ t \in [0, 1] : F( \lambda (t)) \in M \}$ where, as before, $\lambda = \lambda_{a', b'-a'}$ and $a' = \langle a_{i_1} , \ldots , a_{i_m} \rangle$, $b' = \langle b_{i_1} , \ldots , b_{i_m} \rangle$, 
 
 Then $T$ is not empty because $0 \in T$. Also, $T$ is a closed subset of $[0, 1]$ because $M$ is closed (2(a) of Definition 2). So we shall be done if we can show that $T$ is open.
 
 \vspace{2mm}
 
 So let $t_0 \in [0, 1]$ be such that $F(\lambda (t_0 )) \in M$. Say $F(\lambda (t_0 )) = c = \langle c_1 , \ldots , c_n \rangle$. Choose $G = \langle G_1 , \ldots , G_{n-l} \rangle$ and $W$ as in (2b) of Definition 2 for this particular $c \in M$. Then $c \in Z_{reg} (G)$, and by reducing $W$ (definably) if necessary we may suppose that $w$ is a non-singular point of the zero set of $G$ for \emph{all} $w \in W$ satisfying $G(w) =0$.
 By continuity, there exists $\epsilon >0$ such that $F(\lambda (t)) \in W$ for all $t \in [t_0 - \epsilon , t_0 + \epsilon ]$. Thus $\{ z \in V : F(z) \in W \}$ is a definable, open subset of $V$ containing $\lambda ([t_0 - \epsilon , t_0 + \epsilon ])$ and $G \circ F : \{ z \in V : F(z) \in W \}$* $\rightarrow \C$ is a definable holomorphic function such that $G \circ F(\lambda (t_0 )) = 0$, 
 where the * denotes taking the connected component of the set $\{ z \in V : F(z) \in W \}$  that contains the point $\lambda (t_0 )$ (and hence the set $\lambda ([t_0 - \epsilon , t_0 + \epsilon ])$). However, $\lambda$ is generic on $[0, 1]$ and so it follows from (2) that $G \circ F$ is identically zero. In particular $G(F(\lambda (t))) = 0$ for all $t \in [t_0 - \epsilon , t_0 + \epsilon ]$. But $Z_{reg} (G) \subseteq M \cap W$ ((2b) of Definition 2) and so $F(\lambda (t)) \in M$ for all $t \in [t_0 - \epsilon , t_0 + \epsilon ]$, and this shows that $T$ is open, as required.
 
\vspace{3mm}

The proof of our present aim is now complete apart from one small detail. The reader may have noticed that, strictly speaking, atomic formulas of $\mathcal{L} (\widetilde{\C} )$ of $\widetilde{\C}$ have the form $\Phi (v_{j_1} , \ldots , v_{j_p} )$ for some $1 \leq l < p$ where $\Phi$ is the symbol of the language $\mathcal{L} (\C )$ interpreting an $l$-submanifold of $\C^p$ (for some $l, p$ with $1 \leq l < p$). But we have tacitly assumed in our proof above that $j_k = k$ and that $p$ is (an arbitrarily large) $n$. But this assumption can easily be arranged (at least for $j_1 , \ldots , j_p$ distinct) by ``adding vacuous variables'' and observing that the set

\vspace{3mm}

\hspace{15mm} $\{ \langle a_1 , \ldots , a_n \rangle \in \Cn : \tilde{\C} \models \Phi [a_{j_1} , \ldots , a_{j_p} ] \}$

\vspace{3mm}

is an $(n-p+l)$-submanifold of \Cn\ . 

 \vspace{2mm}
 
 Notice also that the graph of equality is a $1$-submanifold of $\C^2$ (so we may indeed assume in the above discussion that $j_1 , \ldots , j_p$ are distinct) and that the graphs of addition and multiplication are $2$-submanifolds of $\C^3$, so that  equality of polynomial terms is also preserved by the $\sim_n$ relation. 
 
 \vspace{2mm}
 
 We now need to deal with case that the formula defining the set $S$ contains a countable set, $X$ say, of parameters. But for this we simply apply the result above with $\widetilde{\R}$ replaced by the structure, $\breve{\R}$ say, obtained by expanding $\widetilde{\R}$ by a constant for each element of $X'$ where $X'$ is the set of real and imaginary parts of elements of $X$. The required result follows since it is easy to check (upon denoting by $\breve{\C}$ the corresponding complex structure as given by Definition 3) that for any formula of $\mathcal{L}(\widetilde{\C})_{\infty , \omega}$ with parameters in $X$ there exists a parameter-free formula of 
 $\mathcal{L}(\breve{\C})_{\infty , \omega}$ that defines the same set. $\Box$
 
 \vspace{4mm}

 As a test for quasi-minimality Theorem 1 has limited use because in order to prove that a given structure has the ACP one still needs some reasonable mathematical description of the definable sets and functions. We now turn to this problem in the case of expansions of the complex field by entire functions and we look for a \emph{complex analytic} criterion for such a structure to have the ACP.
 To this end, suppose that we are given, for each $n \geq 0$, a countable ring $\mathcal{H}_n$ of entire functions of the $n$ complex variable $z_1 , \ldots , z_n$. We assume that $\mathcal{H}_n \subseteq \mathcal{H}_{n+1}$  (in the obvious sense) and that each $\mathcal{H}_n$ contains the projection functions and is closed under partial differentiation and Schwarz reflection (i.e. if $f \in \mathcal{H}_n$ , then $\frac{\partial f}{\partial z_j} \in \mathcal{H}_n$ for $j = 1, \ldots , n$ and $f^{SR} \in   \mathcal{H}_n$, where $f^{SR}(z) := \overline{f(\overline{z} )}$ for $z \in \Cn$ (and the bar denotes co-ordinatewise complex conjugation)).  We then call the sequence $\mathcal{H} := \langle  \mathcal{H}_n : n \geq 0 \rangle$ of rings  a \emph{suitable sequence} and associate to such an  $\mathcal{H}$ a certain expansion $\tilde{\R}( 
 \mathcal{H})$ of $\overline{\R}$ as follows:-
 
 \vspace{3mm}
 
 \noindent
 For each $n \geq 0$, $f \in \mathcal{H}_n$ and discs $D_1 , \ldots , D_n$ in \C\ with Gaussian rational centres and rational radii we denote by $\tilde{f}$ the restriction of $f$ to the polydisc $D_1 \times \cdots \times D_n$ (and define $\tilde{f}(z)$ to be $0$ for $z \notin D_1 \times \cdots \times D_n$). (For $n=0$, $\tilde{f}$ is taken to be the element $f$ of  $\mathcal{H}_0$ ($\subseteq \C$).)

 Now, with the usual convention concerning the identification of \C\ with $\R^2$, we define the structure $\widetilde{\R}( 
 \mathcal{H})$ to be the expansion of $\overline{\R}$ by all such $\tilde{f}$.
 
 \vspace{2mm}
 
 Then $\widetilde{\R}( \mathcal{H})$, being a reduct of $\R_{an}$, is a polynomially bounded, o-minimal expansion of $\overline{\R}$ and its language is countable (since each 
$\mathcal{H}_n$ is and there are only countably many polydiscs to which we restrict the functions therein).

\vspace{2mm}

In [W3] I gave a characterization of the definable, holomorphic functions of  $\widetilde{\R}(\mathcal{H})$ around \emph{generic} points of \Cn\ . This characterization has been shown to be insufficient around non-generic points (see [JKLS]), but at least it does give an alternative  description of the $LD$-pregeometry in terms that avoid notions of general definability. The characterization is as follows.

\vspace{2mm}

Consider an $LD$-generic point $a = \langle a_1 , \ldots , a_n \rangle \in \Cn$ and let $F$ be a definable, holomorphic function (definable, that is, in the structure $\widetilde{\R}( \mathcal{H})$ without parameters). Then, as is proved in [W3], there exist disks $D_1 , \ldots , D_n$ in \C\  (with centres and radii as specified above) and, for some $N \geq 1$, functions $f_1 , \ldots , f_N \in \mathcal{H}_{n+N}$, and definable, holomorphic functions $\phi_1 , \ldots , \phi_N :  D_1 \times \cdots \times D_n \rightarrow \C$ such that

\vspace{3mm}

\noindent  
(3) \hspace{2mm} $a \in D_1 \times \cdots \times D_n \subseteq dom(F)$;
 
\vspace{3mm}

\noindent
(4) \hspace{2mm} for all $z \in D_1 \times \cdots \times D_n$ we have that
 $\langle \phi_1 (z), \ldots , \phi_N (z) \rangle$ is a non-singular zero of the map $f_z : \C^N \rightarrow \C^N : w \mapsto \langle f_1 (z, w), \ldots , f_N (z, w) \rangle$;
 
 \vspace{3mm}
 
 \noindent
 (5) \hspace{2mm} $F(z) = \phi_1 (z)$ for all $z \in   D_1 \times \cdots \times D_n$.
 
\vspace{3mm}
 
 In other words, $F$ arises, at least close to the generic point $a$, as a coordinate function of a map given by an application of the Implicit Function Theorem applied to functions from $\mathcal{H}$. (Actually, the characterization from [W3] (see 1.5 and 1.6 of [W3]) makes use of just the one (dependent) variable version of the Implicit Function Theorem together with composition, but it is easy to see that this formulation is equivalent. Note also that the operations 1.2 and 1.3 from [W3] follow from the corresponding closure conditions that we have placed on $\mathcal{H}$. (One can consult [S] for more on this,))
 
\vspace{3mm}
 
 For later use we remark now that (4) is equivalent to
 
 \vspace{3mm}

 \noindent  
 (4*) \hspace{2mm} for all $z \in D_1 \times \cdots \times D_n$ and $j = 1, \ldots , N$ we have that \\ $f_j (z , \phi_1 (z) , \ldots , \phi_N (z) ) = 0$ and $J(z, \phi_1 (z) , \ldots , \phi_N (z)) \neq 0$ where $J$ is the determinant of the Jacobian matrix $( \frac{\partial f_j}{\partial w_i} )_{1 \leq i, j \leq N}$. (Note that $J \in \mathcal{H}_{n+N}$.)
 
 \vspace{4mm}
 
 In [W3] I define a pregeometry $\widetilde{D}$ on \C\  associated with implicit functions as discussed above. Namely, a $d$-tuple $\langle b_1 , \ldots , b_d \rangle \in \C^d$   is declared to be $\widetilde{D}$-generic if there do not exist $k \geq 0$ and $b_{d+1} , \ldots , b_{d+k} \in \C$ and 
 $g_{1} , \ldots , g_{k+1} \in \mathcal{H}_{d+k}$ such that $\langle b_1 , \ldots , b_d , b_{d+1} , \ldots , b_{d+k} \rangle$ is a non-singular zero of the map $\langle g_{1} , \ldots , g_{k+1} \rangle : \C^{d+k} \rightarrow \C^{k+1}$.
 
 It is shown in [W3] (Theorem 1.10) that $LD$ and $\widetilde{D}$ are identical pregeometries and as regards to our present aim this leads to our next theorem. It states that one only needs to check that implicitly defined functions, such as the $\phi_i$'s mentioned above, have 
 analytic continuations along generic paths. To be more precise we make the following
 
 \vspace{3mm}
 
 \noindent
 $\mathbf{Definition}$ $\mathbf{4}$
 
 \vspace{1mm}
 
 We say that the suitable sequence $\mathcal{H}$ has the \emph{Weak Analytic Continuation Property} (WACP) if for all $n$, $N \geq 1$, all $a$, $\omega \in \Cn$ such that $\lambda_{a, \omega}$ is $\widetilde{D}$-generic on $[0, 1]$, all $r \in [0, 1]$ and all $f \in (\mathcal{H}_{n+N})^N$, if $\gamma : [0, r) \rightarrow \C^N$ is a continuous map such that $\langle \lambda_{a, \omega} (t), \gamma (t) \rangle \in Z_{reg} (f)$ for all $t \in [0, r)$ then $|| \gamma (t) || \nrightarrow \infty$ as $t \rightarrow \infty$. (The non-singularity here is with respect to the last $N$ variables (as in (4*)), and $|| \cdot ||$ is some standard norm on $\C^N$.)
 
 \vspace{3mm}
 
 Notice that there is no explicit mention of definability here. Nevertheless, we have the following
 
 \vspace{3mm}
 
 \noindent 
 $\mathbf{Theorem}$ $\mathbf{2}$
 
 \vspace{1mm}
 
 If $\mathcal{H}$ has the WACP then $\widetilde{\R}(\mathcal{H})$ has the ACP. Hence, by Theorem 1, the corresponding expansion $\widetilde{\C}(\mathcal{H})$ of the complex field is 
 quasiminimal (even for the language $\mathcal{L}(\widetilde{\C}(\mathcal{H})_{\infty , \omega}$).
 
 \vspace{3mm}
 
 \noindent
$\mathbf{Reminder}$ $\mathbf{and}$ $\mathbf{Remark}$
 
\vspace{1mm}

The structure $\widetilde{\C}(\mathcal{H})$
 is defined to be the expansion of the complex field by all $l$-dimensional, locally definable, complex submanifolds of \C\  (for all $l$, $n$ with $1 \leq l < n$) with respect to the expansion $\widetilde{\R}(\mathcal{H})$ of the real field. It is trivial to show that all functions in $\mathcal{H}$ are definable in $\widetilde{\C}(\mathcal{H})$.

\vspace{3mm}

 \noindent
 $\mathbf{Proof}$ $\mathbf{of}$ $\mathbf{Theorem}$ $\mathbf{2}$
 
 \vspace{1mm}
 
 Let $a \in \Cn$ be an $LD$-generic point and $F$ a definable holomorphic function with $a \in dom(F)$. Let $\omega \in \Cn$ be such that $\lambda_{a, \omega}$ is generic on $[0, 1]$. We write $\lambda$ for 
 $\lambda_{a, \omega}$. We must find $G$ satisfying the conclusion of Definition 1.
 
 Choose $\langle D_1 , \ldots , D_n \rangle$, $N$, $f = \langle f_1 , \ldots , f_N \rangle$ and $\phi = \langle \phi_1 , \ldots , \phi_N \rangle$ as in (3), (4) and (5).
 
 Now we may suppose that \emph{every} zero, $\langle z^{(0)}, w^{(0)} \rangle$ say, of $f$ satisfies \\ $J(\langle z^{(0)}, w^{(0)} \rangle ) \neq 0$ (see (4*)). Indeed, let $f_{N+1} \in \mathcal{H}_{n+N+1}$ be defined by
 
 \vspace{3mm} 
 
 \hspace{25mm} $f_{N+1}(z, w, w_{N+1}) := w_{N+1} \cdot J(z, w) - 1$.
 
 \vspace{3mm}
 
 \noindent
 Then, letting $f^* := \langle f_1 , \ldots , f_N , f_{N+1} \rangle$, one easily calculates that the Jacobian of $f^*$ (with respect to $w_1 , \ldots , w_{N+1}$) has determinant $J(z, w)^2$, which is non-zero whenever $f_{N+1} (z, w, w_{N+1})$ is zero. Further, any non-singular (with respect to $w_1 , \ldots , w_{N+1}$) zero of $f$ gives rise to a (unique) zero of the map $f^*$. So we may replace $f$ by $f^*$ and (3), (4) and (5) remain true (by setting $\phi_{N+1}(z)  = J(z, \phi_1 (z), \ldots , \phi_N (z))^{-1}$). And now, \emph{all} zeros of $f^*$ are non-singular (with respect to $w_1 , \ldots , w_{N+1}$).
 
 So we continue our proof with this non-singularity assumption.
 
 \vspace{2mm}
 
 Let $T$ be the set of all those $t \geq 0$ having the following properties:-
 
 \vspace{2mm}
 
 \noindent
$(6)_t$ \hspace{2mm} there exists a definable, open, connected set $U_{t} \subseteq \C$ with $[0, t] \subseteq U_t$;
 
 \vspace{2mm}
 
 \noindent
$(7)_t$ \hspace{2mm} there exists a definable, holomorphic map $\psi^{(t)}$ with range contained in $\C^N$ and with $\lambda (U_t ) \subseteq dom(\psi^{(t)})$ which satisfies
 
 \vspace{2mm}
 
\noindent
$(8)_t$ \hspace{2mm} for all $u \in U_t$ we have that $\langle \lambda (u) , \psi^{(t)} (\lambda (u) \rangle$ is a non-singular zero of $f$ (with respect to $w_1 , \ldots , w_{N}$);

\vspace{2mm}

\noindent
$(9)_t$ \hspace{2mm} $\psi^{(t)} (\lambda (0)) = f( \lambda (0))$. (So, in particular, the first coordinate of $\psi^{(t)} (\lambda (0))$ is $F(a)$.)
 
 \vspace{3mm}
 
 We shall be done if we can show that $1 \in T$, for then we take $G$ to be the first coordinate of $\psi^{(1)}$ to satisfy the conclusion of Definition 1.
 
 \vspace{3mm}
 
 Notice that we certainly have $0 \in T$ by taking $U_0$ to be a (definable) disk around $0 \in \C$ which is small enough to satisfy $\lambda (U_0 ) \subseteq D_1 \times \cdots \times D_n$, and then taking $\psi^{(0)} := \phi$.
 
 Notice also that if $t_1$, $t_2 \in T$ and $t_1 \leq t_2$ then $\psi^{(t_1 )} \circ \lambda$ and $\psi^{(t_2 )} \circ \lambda$ must agree on $(U_{t_1} \cap U_{t_2})^*$, the connected component of $U_{t_1} \cap U_{t_2}$ containing the interval $[0, t_1 ]$. This is because they agree at $0$ (by $(9)_{t_1}$ and $(9)_{t_2}$) and $\lambda (0)$ ($= a$) is $LD$-generic, so the definable holomorphic map \\$(\psi^{(t_2 )} - \psi^{(t_1 )})$$\upharpoonright$$(U_{t_1} \cap U_{t_2})^*$ must be identically zero (the zero, that is, of $\C^N$).
 
 \vspace{1mm}
 
 We now set $r := sup \{ t : [0, t] \subseteq T \}$ and we need to show that $r \geq 1$. So suppose, for a contradiction, that $r < 1$. By the extension property just proved, it follows that $\chi := \bigcup_{t < r} (\psi^{(t)} \upharpoonright \lambda ([0, t]))$ is a continuous map with domain $\lambda ([0, r))$ such that for all $t \in [0, r)$ we have that $\langle \lambda (t), \chi (\lambda (t)) \rangle$ is a non-singular zero of $f$ (with respect to $w_1 , \ldots , w_N$) - see $(8)_t$. So by applying the WACP (with $\gamma = \chi \circ \lambda$) we see that there exists some positive $R$ and an increasing sequence $\langle t_p : p \geq 0 \rangle$ in $[0, r)$ converging to $r$ such that $|| \psi ( \lambda (t_p ))|| \leq R$ for all $p \geq 0$.
 
 Let $w^{(0)} \in \C^N$ be a limit point of the sequence $\langle \chi ( \lambda (t_p )) : p \geq 0 \rangle$. Then $ \langle \lambda (r) , w^{(0)} \rangle$ is a zero of $f$ (since $f$ is certainly a continuous map throughout $\C^{n+N}$) and by our non-singularity assumption, it is a non-singular zero (with respect to $w_1 , \ldots , w_N$). So by the Implicit Function Theorem there exist an open polydisc $\Delta \subseteq \Cn$ (which we may take to be definable) with $\lambda (r) \in \Delta$, and a holomorphic map $\theta : \Delta \rightarrow \C^N$ satisfying $\theta (\lambda (r)) = w^{(0)}$ and such that for all $z \in \Delta$, the $(n+N)$-tuple $\langle z, \theta (z) \rangle$ is a non-singular zero of $f$ (with respect to $w_1 , \ldots , w_N$). Further, we may assume that $\Delta$ has been chosen small enough for there to exist a (definable) open polydisk $E \subseteq \C^N$ with $w^{(0)} \in E$ such that for all $z \in \Delta$, $w = \theta (z)$ is the one and only solution in $E$ of the equation $f(z, w) = 0$. (This follows from the uniqueness condition in the conclusion of the Implicit Function Theorem.) It follows from this that $\theta$ is definable.
 
 \vspace{1mm}
 
 Now choose $p$ large enough so that $t_p$ is close enough to $r$ to satisfy 
 
 \vspace{3mm}
 
 \noindent
 (10) \hspace{2mm} $\lambda (t_p ) \in \Delta$ and
 
 \vspace{3mm}
 
 \noindent
 (11) \hspace{2mm} $\chi (\lambda (t_p )) \in E$ (and hence $\psi^{(t_p )} (\lambda (t_p )) \in E$).
 
 \vspace{3mm}
 
 Now choose $\epsilon > 0$ so small that the rectangle
 
 \vspace{3mm}
 
 \hspace{20mm} $\rho := \{ x + iy \in \C : t_p -\epsilon < x < t_p + \epsilon , -\epsilon < y < \epsilon \}$
 
 \vspace{2mm}
 
 \noindent
 is  contained in $U_{t_p}$. (See $(6)_{t_p}$.)
 
 \vspace{2mm}
 
 Since $\lambda$ is linear, $\lambda (\rho )$ is a convex, open subset of \Cn  (=$\R^{2n}$) and since $\Delta$ is too, it follows that $\lambda (\rho ) \cap \Delta$ is convex and, in particular, connected.
 
 Now $\lambda (t_p ) \in \lambda (\rho ) \cap \Delta$ (see (10)) and both $\theta (\lambda (t_p ))$ and $\phi^{(t_p )} (\lambda (t_p ))$ lie in $E$ (see (11)). So by the uniqueness condition we have that $\phi^{(t_p )} (\lambda (t_p )) = \theta (\lambda (t_p 
 ))$ (see $(8)_{t_p}$). 
 
 But $\lambda$ is generic (for either pregeometry) on $[0, 1]$ so, by (2), $\phi^{(t_p )}$ and $\theta$ must agree on a sufficiently small, definable open polydisk containing the point $\lambda (t_p )$ and contained within $\lambda (\rho ) \cap \Delta$. But then, by the principle of analytic continuation, they must agree throughout the connected set 
 $\lambda (\rho ) \cap \Delta$. (Note that $\lambda (\rho ) \cap \Delta \subseteq dom(\phi ^{(t_p)}) \cap dom (\theta )$-see $(7)_{t_p}$.)
 So we may consistently define a holomorphic map, $\Gamma$ say,  with domain 
  $\lambda (\rho ) \cup \Delta$ and taking values in $C^N$ by specifying $\Gamma (z)$ to be $\phi^{(t_p )} (z)$ for $z \in \lambda (\rho )$ and $\theta (z)$ for $z \in \Delta$.
  
  Finally, since $\lambda (r) \in \Delta$ we may choose $r'$ with $r < r' \leq 1$ such that $\lambda (r') \in \Delta$. We now obtain a contradicion by showing that $[0, r'] \subseteq T$. 
  
  Indeed, suppose that $s \in [0, r']$. Then $\lambda ([0, s]) \subseteq \lambda (\rho ) \cup \Delta$ (because $\lambda ([0, t_p ]) \subseteq \lambda (\rho )$ and $\lambda ([t_p ,s]) \subseteq \Delta$ (since $\Delta$ is convex)). So if we take $U_s$ to be a definable, open, connected subset of \C\  containing the interval $[0, s]$ and contained within $\lambda^{-1} (\lambda (\rho ) \cup \Delta )$, and set $\psi^{(s)} := \Gamma \upharpoonright \lambda (U_s )$, we see that $(6)_s$-$(9)_s$ are all satisfied. So $s \in T$ as required. $\Box$ 
 
 \vspace{5mm}
 
 \noindent
 $\mathbf{Exponentiation}$
 
 \vspace{2mm}
 
 We now consider the setting appropriate, with respect to the preceeding discussion of analytic continuation, for Zilber's quasiminimality conjecture for the complex exponential field.
 
 \vspace{1mm}
 
 We let K be a fixed countable, real closed subfield of \R\  and let $\mathcal{E}_0$ be the algebraically closed subfield $K[i]$ of \C\ . Then, for each $n \geq 1$ we let $\mathcal{E}_n$ be the ring $\mathcal{E}_0 [ z_1 ,  \ldots , z_n , e^{z_1} , \ldots , e^{z_n} ]$ of \emph{exponential polynomials} over $\mathcal{E}_0$ in the complex variables $z_1 ,  \ldots , z_n$. It is clear that $\mathcal{E} := \langle 
 \mathcal{E}_n : n \geq 0 \rangle$ is a suitable sequence and one easily checks that $\widetilde{\R} ( \mathcal{E} )$ is essentially (i.e. has the same definable sets as) as the structure   $\R^{RE}_{K}$ considered by Binyamini and Novikov in [BN]. The superscript  ``$RE$'' stands for ``restricted elementary'':  $\R^{RE}_{K}$ is the expansion of $\overline{\R}$ by  the restricted functions $\exp$$\upharpoonright$$[0, 1]$ and $\sin$$\upharpoonright$$[0, \pi ]$ and by a constant for each element of  $K$.

 \vspace{1mm}
 
 Thus, in the present context, a subset $M$ of \Cn\  is an $l$-dimensional, locally definable, complex submanifold of \Cn\  if it satisfies (2a) and (2b) of Definition 2, where definability is now with respect to the structure  $\R^{RE}_{K}$ (and, as before, is \emph{without} parameters).
 
 \vspace{1mm}
 
 Let us refer to such $M$ simply as \emph{elementary complex manifolds} (ECM's). We write $\C^{ECM}_{K}$ for the structure denoted $\widetilde{\C}$ in Definition 3; that is, $\C^{ECM}_{K}$ is the expansion of the complex field by all ECM's.
 
 $\C^{ECM}_{K}$ is certainly an expansion of $\C_{exp}$. (See the remark following Definition 3.)
 Notice also that the notion of an ECM is, at least apparently, more general than if we had required the functions $G_1 , \ldots , G_{n-l}$ of (2b) (of Definition 2) to lie in $\mathcal{E}_n$. Perhaps the methods of [JKLS] could be used to produce an ECM which is not of this latter kind.
 
 Theorem 2 implies
 
 \vspace{2mm}
 
 \noindent
 $\mathbf{Theorem}$ $\mathbf{3}$
 
 \vspace{1mm}
 
 If $\mathcal{E}$ has the WACP then the structure $\C^{ECM}_{K}$ (and so, in particular, the structure $\C_{exp}$) is quasiminimal (even for the language $\mathcal{L}(\C^{ECM}_{K})_{\infty , \omega}$).
 
 \hspace{120mm} $\Box$
 
 \vspace{4mm}
 
 I conclude this article with a proof that $\mathcal{E}$ has the WACP in the case $N=1$ in the hope that others might find a way to generalize the method, which uses the (algebraic) valuation inequality.
 
 \vspace{2mm}
 
 So, suppose that $n \geq 1$ and $a, \omega \in \Cn$ are such that $\lambda$ := $\lambda_{a, \omega}$ is $\widetilde{D}$-generic on $[0, 1]$. Suppose further that $f \in \mathcal{E}_{n+1}$, $r \in [0, 1]$ and that $\gamma : [0, r) \rightarrow \C$ is a continuous function such that
 
 \vspace{3mm}
 
 \noindent
 (12) \hspace{2mm} $f(\lambda (t), \gamma (t)) = 0 \neq \frac{\partial f}{\partial z_{n+1}} (\lambda (t), \gamma (t))$
 
 \vspace{1mm}
 
 \noindent
 for all $t \in [0, r)$.
 
 \vspace{2mm}
 
 We must show that $|\gamma (t)| \nrightarrow \infty$ as $t \rightarrow r$.
 
 \vspace{1mm}
 
 In order to set up a use of the valuation inequality we require the following general
 
 \vspace{2mm}
 
 \noindent
 $\mathbf{Lemma}$
 
 \vspace{1mm}
 
 Suppose that $r > 0$ and that $\phi : [0, r) \rightarrow \C$ is a continuous function such that $|\phi (t)| \rightarrow \infty$ as $t \rightarrow r$. Then \emph{either} (A) for all integers $k$, $l$, not both $0$ with $k \geq 0$, we have either $|\phi (t)^k \exp(l\phi (t))| \rightarrow 0$ as $t \rightarrow r$ or $|\phi (t)^k \exp(l\phi (t))| \rightarrow \infty$ as $t \rightarrow r$,  \emph{or} (B) for some integers 
$k$, $l$, not both $0$ with $k \geq 0$ we have that for all countable sets $S \subseteq \C$, there exists $\alpha \in \C \setminus S$ and an increasing sequence $0 \leq t_0 < t_1 < \ldots < t_p < \ldots$ converging to $r$ such that $\lim_{j \rightarrow \infty} \phi (t_j)^k \exp(l\phi (t_j)) = \alpha$.
 
\vspace{1mm}

\noindent
$\mathbf{Proof}$

\vspace{1mm}

Set $J := \{ \langle k, l \rangle \in \mathbb{Z}^2 : k \geq 0$ and $k$, $l$ not both $0 \}$.

For $\langle k, l \rangle \in J$ write $h_{k, l} (t)$ for  $\phi (t)^k \exp(l\phi (t))$ (for $t \in [0, r)$), and define $c_{k, l}^{+} := \limsup_{t \rightarrow r}  |h_{k, l} (t)|$  and $c_{k, l}^{-} := \liminf_{t \rightarrow r}  |h_{k, l} (t)|$.  

Then $0 \leq  c_{k, l}^{-} \leq c_{k, l}^{+} \leq \infty$.

If for each $\langle k, l \rangle \in J$ we have either $0 = c_{k, l}^{-} = c_{k, l}^{+}$ or $c_{k, l}^{-} = c_{k, l}^{+} = \infty$ then clearly (A) holds.

Otherwise, choose $\langle k, l \rangle \in J$ with $c_{k, l}^{+} >0$ and $c_{k, l}^{-} < \infty$. Let $S$ be a countable subset of $\C$. Write $h$ for $h_{k, l}$.

\vspace{1mm}

Now either $c_{k, l}^{-} < c_{k, l}^{+}$ or $0 < c_{k, l}^{-} = c_{k, l}^{+} < \infty$.

In the first case choose $c \in \R$ with $c_{k, l}^{-} < c < c_{k, l}^{+}$ and $c \notin \{ |s| : s \in S \}$.

 By the continuity of $|h|$ there clearly exists a sequence $0 \leq t'_0 < t'_1 < \ldots$ converging to $r$ such that $|h(t'_j)| = c$ for all $j \in \mathbb{N}$. But now \\ $\langle h(t'_j ) : j \in \mathbb{N} \rangle$ is a bounded sequence of complex numbers and hence has a convergent subsequence $\langle h(t_j ) : j \in \mathbb{N} \rangle$ whose limit, $\alpha$ say, cannot lie in $S$ because $|\alpha | = c$.

\vspace{2mm}

In the second case we have $\lim_{t \rightarrow r} |h(t)| = c$, say, with $0 < c < \infty$.
There is no harm in assuming that both $h$ and $\phi$ are nonzero throughout $[0, r)$ and so there exist continuous functions $\theta$, $\psi : [0, r) \rightarrow \R$ such that

\vspace{2mm}

\hspace{25mm} $h(t) = |h(t)|\exp (i\theta (t))$, and

\vspace{1mm}

\hspace{25mm} $\phi (t) = |\phi (t)|\exp(i\psi (t))$,

\vspace{1mm}

\noindent
for all $t \in [0, r)$.

\vspace{2mm}

So by definition of $h$ we have 

\vspace{2mm}

\noindent
(13) \hspace{3mm} $|h(t)|\exp(i\theta(t)) = |\phi (t)|^k \exp (ik\psi (t)) \cdot \exp (l|\phi (t)|(\cos \psi (t) + i\sin \psi (t)))$. 

\vspace{2mm}

Hence $|h(t)| = |\phi (t)|^k \exp (l|\phi (t)|\cos( \psi (t)) \rightarrow c$ as $t \rightarrow r$.

\vspace{1mm}

We cannot have $l = 0$, for then $k > 0$ and $\phi$ would be bounded. So $l \neq 0$ from which it follows that $\cos \psi (t) \rightarrow 0$ as $t \rightarrow r$ (since $c \neq 0, \infty$). Thus $\psi$ is bounded and $\sin \psi (t) \rightarrow \pm 1$ as $t \rightarrow r$. Equating arguments in (13) we obtain, for some fixed $N_0 \in \mathbb{Z}$ and for all $t \in [0, r)$:

\vspace{3mm}

\hspace{15mm} $\theta (t) = k\psi (t) + l|\phi (t)|\sin \psi (t) + 2\pi N_0$.

\vspace{3mm}

It follows from this that $\theta (t) \rightarrow \pm \infty$ (the sign here depending on the eventual sign of $l\sin \psi (t)$) as $t \rightarrow r$.

Thus we may choose some $\theta_0 \in \R$ such that $c \exp (i \theta_0 ) \notin S$ and for which there exists a sequence $0 \leq t_0 < t_1 < \ldots$ converging to $r$ such that $\theta (t_j) = \theta_0$ (mod $2\pi \mathbb{Z}$) for all $j \in \mathbb{N}$. It follows that $h(t_j ) \rightarrow c\exp (i\theta_0 )$ as $j \rightarrow \infty$, and we are done. $\Box$.

\vspace{3mm}

Now returning to the discussion before the statement of the lemma, suppose, for a contradiction, that (12) holds and that $|\gamma (t)| \rightarrow \infty$ as $t \rightarrow r$. By definition of $\mathcal{E}_{n+1}$ we see that $f$ has the form

\vspace{3mm}

\noindent
(14) \hspace{5mm} $f(z_1 , \ldots , z_n , z_{n+1} ) = \sum_{\langle i, j \rangle \in L} P_{i, j} (z_1 , \ldots , z_n ) z_{n+1}^i \exp (jz_{n+1})$

\vspace{3mm}

\noindent
for some non-empty finite set $L \subseteq \mathbb{N}^2$, where $P_{i, j} \in \mathcal{E}_n \setminus \{ 0 \}$ for each $\langle i, j \rangle \in L$. We must have $L \neq \{ \langle 0, 0 \rangle \}$ by (12).

By the genericity of $\lambda$ on $[0, r]$ it routinely follows that for all $P \in \mathcal{E}_n \setminus \{ 0 \}$ there exists some $R_P \geq 1$ such that

\vspace{3mm}

\noindent
(15) \hspace{10mm}  $R_P \geq |P(\lambda (t))| \geq R^{-1}_P$ for all $t \in [0, r]$.

\vspace{3mm}

Let us pass to a non-principal ultrapower $^{\star} \C$ of \C\  (with corresponding
 $^{\star}\R$, $^{\star}\mathbb{Z}$, $^{\star} \mathbb{N}$). Then the functions $\lambda_1 , \ldots , \lambda_n$ (where $\lambda = \langle \lambda_1 , \ldots , \lambda_n \rangle$) and all the $P_{i, j}$'s have natural extensions to the ultrapower and (keeping the same notation for the extended functions) (15) remains true for all $t \in$ $^{\star}\R$ with $0 \leq t \leq r$.
 
 For each such $t$ consider the subfield
 
 \vspace{2mm}
 
 \hspace {4mm} $\mathcal{F}_t  := \mathcal{E}_0 ( \lambda_1 (t) , \ldots , \lambda_n (t), \exp (\lambda_1 (t)) , \ldots , \exp(\lambda_n (t)))$
 
 \vspace{2mm} 
 
 \noindent
 of $^{\star} \C$.
 
 \vspace{2mm}
 
 Then by (the extension to the ultrapower of) (15) it follows that $\mathcal{F}_t$ is actually a subfield of the valuation subring Fin($^{\star}\C$) (:= $\{ z \in ^{\star}\C : |z| \leq R$ for some $R \in \R \}$) of $^{\star}\C$.

 By the continuity of each $P \circ \lambda$ (for $P \in \mathcal{E}_n$) it follows that for all $t_1$, $t_2 \in ^{\star}\R$ with $0 \leq t_1 , t_2 \leq r$ and satisfying $t_1 \approx t_2$ (i.e. $t_1$ infinitesimally close to $t_2$) we have that $P( \lambda (t_1 )) \approx P(\lambda (t_2 ))$ and so the correspondence $t_1 \mapsto t_2$ induces an isomorphism $\mathcal{I}_{t_1 ,t_2} : \mathcal{F}_{t_1} \rightarrow \mathcal{F}_{t_2}$ with $\mathcal{I}_{t_1 ,t_2} (z) \approx z$ for all $z \in \mathcal{F}_{t_1}$ (and so, in particular, $\mathcal{I}_{t_1 ,t_2} (z) = z$ for all $z \in \mathcal{E}_0$).

 Further, by the continuity of roots of polynomials (see [HM] or, perhaps more appropriately in the present context, [R]) the map 
 $\mathcal{I}_{t_1 ,t_2}$ extends to an isomorphism $\widetilde{\mathcal{I}}_{t_1 ,t_2} : \widetilde{\mathcal{F}}_{t_1} \rightarrow \widetilde{\mathcal{F}}_{t_2}$ of the algebraic closures and we still have (for $t_1 \approx t_2$):
 
 \vspace{2mm}
 
 \noindent
 (16) \hspace{2mm} (a) $\widetilde{\mathcal{F}}_{t_1}$, $\widetilde{\mathcal{F}}_{t_2} \subseteq$ Fin$(^* \C$) and (b) $\widetilde{\mathcal{I}}_{t_1 ,t_2} (z) \approx z$ for all $z \in \widetilde{\mathcal{F}}_{t_1}$.
 
 \vspace{2mm}
 
 Now choose any $t^* \in$ $^*\mathbb{R}$ with $0 < t^* < r$ and $t^* \approx r$. Extend the function $\gamma : [0, r) \rightarrow \C$ to the ultrapower and set
  
 \vspace{2mm}
 
 \noindent
 (17) \hspace{20mm} $\mathcal{H}_{t^*} := \widetilde{F}_{t^*} (\gamma (t^* ), \exp (\gamma (t^* )))$.
 
 \vspace{2mm}
 
 Then by (12) and (14), $\mathcal{H}_{t^*}$ is a subfield of $^*\C$ of transcendence degree at most $1$ over $\widetilde{\mathcal{F}}_{t^*}$ and, in fact, exactly $1$ because $\gamma (t^* ) \notin$ Fin($^*\C$) (and $\widetilde{\mathcal{F}}_{t^*} \subseteq$ Fin($^*\C$)). It also follows from this (and the valuation inequality for the valuation on $^*\C$ (restricted to $\mathcal{H}_{t^*}$) associated to the valuation subring Fin($^*\C$) of $^*\C$) that $\gamma (t^* )$ and $\exp (\gamma (t^* ))$ have $\mathbb{Q}$-dependent valuations which, by an easy saturation argument, implies that (A) of the Lemma (with $\phi = \gamma$) cannot hold (back in the standard situation). So (B) holds and we choose integers $k$, $l$ not both $0$ with $k \geq 0$ and take $S$ to be the subset $\widetilde{F}_r$ of \C\ . Let $\alpha \in \C \setminus S$ and $0 \leq t_0 < t_1 < \ldots$ be a sequence converging to $r$ such that $\lim_{j \rightarrow \infty} \phi (t_j )^k \exp (l \phi (t_j )) = \alpha$.

 \vspace{1mm}
 
 Let $p \in$ $^*\mathbb{N} \setminus \mathbb{N}$ and set $t^* := t_p$.
 Then $t^* \approx r$, $t^* < r$ and $\alpha^* := \gamma (t^* )^k \exp (l\gamma (t^* )) \approx \alpha$.
 
 Now for no $\beta \in \widetilde{\mathcal{F}}_{t^*}$ do we have $\alpha^* \approx \beta$; for this would imply $\widetilde{\mathcal{I}}_{t^* , r} (\beta ) \approx \beta \approx \alpha^* \approx \alpha \in \C \setminus \widetilde{\mathcal{F}}_r$, which is absurd since $\widetilde{\mathcal{I}}_{t^* , r} (\widetilde{\mathcal{F}}_{t^*} ) = \widetilde{\mathcal{F}}_r$.
 
 \vspace{1mm}
 
 Since $\widetilde{\mathcal{F}}_{t^*}$ is an algebraically closed subfield of $^*\C$ contained in Fin($^*\C$) this now implies that $\widetilde{\mathcal{F}}_{t^*}(\alpha^{*})$ is also a subfield of $^*\C$ contained in Fin($^*\C$) and is of transcendence degree $1$ over 
 $\widetilde{\mathcal{F}}_{t^*}$. However, certainly $\widetilde{\mathcal{F}}_{t^*}(\alpha{^*}) \subseteq \mathcal{H}_{t^*}$. Further, $\gamma (t^* ) \in \mathcal{H}_{t^*}$ and $\gamma (t^* )$ , being 
infinite, is transcendental over 
 $\widetilde{\mathcal{F}}_{t^*}(\alpha{^*})$ which forces $\mathcal{H}_{t^*}$ to have transcendence degree at least $2$ over $\widetilde{\mathcal{F}}_{t^*}$. This contradiction shows that it cannot be the case that  $|\gamma (t)| \rightarrow \infty$ as $t \rightarrow r$ and completes the proof of the WACP in this rather simple situation. 
 
 \vspace{4mm}
 
 With a little more care one can show that the function $\gamma$ above has an extension to a definable (in the structure $\R_{K}^{RE}$), holomorphic function with $[0, 1] \subseteq$ dom$(\gamma )$ which (therefore) satisfies (12) for all $t \in$ dom($\gamma$). It is then an easy matter to set up an inductive argument (using the full version of the valuation inequality for polynomially bounded o-minimal structures, see for example [Sp]) to establish the WACP for diagonal systems of exponential polynomial equations:
 
 \vspace{2mm}
 
 \noindent
 $\mathbf{Proposition}$

 $\mathcal{E}$ does have the \emph{diagonal} WACP, i.e. where the map $f = \langle f_1 , \ldots , f_N \rangle$ in Definition 4 satisfies the extra condition that $f_j \in \mathcal{E}_{n+j}$ for $j = 1, \ldots , N$. Further, the map $\gamma$ has an extension to a definable (in $\R^{RE}_{K}$), holomorphic function with $[0, 1] \subseteq dom(\gamma)$. $\Box$
 
 \vspace{3mm}
 
 But, unfortunately, this is a \emph{very} special case.
 
 \vspace{5mm}
 
 \noindent
 $\mathbf{Further}$ $\mathbf{remarks}$ $\mathbf{on}$ $\mathbf{quasiminimality}$
 
 \vspace{2mm}
 
 In his paper [B], Boxall shows that every formula (parameters allowed) of the language $\mathcal{L}(\C_{exp})$  having the form $\exists \bar{z} ( P( \bar{w}, \bar{z} ) = 0)$, where $P(\bar{w}, \bar{z})$ is a term of this language (the $\bar{w},  \bar{z}$ being sequences of variables, not necessarily of the same length), is equivalent (in $\C_{exp}$) to a countable boolean combination of formulas of the form 
 $(\exists \bar{z} \in \mathbb{Q}^{m})  \phi(\bar{w}, \bar{z})$ where $\phi$ is a quantifier-free formula of $\mathcal{L}(\C_{exp})$ (containing no parameters other than those used in $P$). This immediately implies that sets of the form 
 
 \vspace{2mm}
 
 \hspace{16mm} $\pi_1(Z(P)) := \{ w \in \C : \exists \bar{z} \in \C^m  P( w, \bar{z}) = 0 \}$
 
 \vspace{2mm}
 
 \noindent
 are either countable or cocountable. 
 
 It is worth mentioning here that, even for the case $m=1$, this is not a property of entire functions in general. For in the paper [A], Alexander complements  earlier work of Tsuji  ([T]) by giving a complete characterization of sets of the form $\pi_1(Z(F))$ for $F: \C^2 \mapsto \C$ an entire function and this characterization implies that there do exist such $F$ with both  $\pi_1(Z(F))$ and $\C \setminus \pi_1(Z(F))$ uncountable. In particular, there is an expansion $\langle \overline{\C}, F \rangle$ of the complex field  $\overline{\C}$ by an entire function $F$ of two variables which is not quasiminimal. However, as pointed out by P. Koiran, it is still not known whether there exists a non-quasiminimal such expansion by an entire function of one variable, or even by finitely many entire functions of one variable. On the other hand, at least we do know (using a combination of ideas from [Ko], [W2] and [Z]) that there exists a \emph{transcendental} entire function $f:\C \mapsto \C$ such that the expansion $\langle \overline{\C}, f \rangle$ of $\overline{\C}$ \emph{is} quasiminimal.

 \vspace{5mm}
 
 \hspace{50mm} $\mathbf{References}$
 
 \vspace{1mm}
 
 \noindent
 [A] H. Alexander, 'On a problem of Julia', Duke Mathematical Journal, Vol. 42, No. 2, (1975), 327-332.
 
 \vspace{1mm}

 \noindent
 [B] Ricardo Bianconi, 'Undefinability results in o-minimal expansions of the real numbers', Annals of Pure and Applied Logic, Vol. 134, Issue 1, June 2005, 43-51.
 
 \vspace{1mm}
 
 \noindent
 [Bo] Gareth Boxall, 'A special case of quasiminimality for the complex exponential field', The Quarterly Journal of Mathematics, 71(3), (2020), 1065-1068.

 \vspace{1mm}

 \noindent
 [BKW] Martin Bays, Jonathan Kirby, A.J. Wilkie, 'A Schanuel property for exponentially transcendental powers', Bull. Lond. Math. Soc. 42 (2010), no.5, 917-922.
 
\vspace{1mm}
 
\noindent
 [BN] Gal Binyamini and Dmitry Novikov, 'Wilkie's conjecture for restricted elementary functions', Annals of Mathematics Vol. 186 (2017), Issue 1, 237-275.
 
\vspace{1mm}
 
\noindent
 [DD] J. Denef and L. van den Dries, $p$-adic and real subanalytic sets, Annals of Mathematics 128 (1988), 79-138.

\vspace{1mm}
 
\noindent
 [HM] Gary Harris and Clyde Martin, 'The roots of a polynomial vary continuously as a function of the coefficients', Proc. Amer. Math. Soc. Vol 100, Number 2, (1987), 390-392.
 
 \vspace{1mm}
 
 \noindent
 [JKLS] Gareth Jones, Jonathan Kirby, Olivier Le Gal and Tamara Servi, 'On local definability of holomorphic functions', Quarterly Journal of Mathematics, 70(4), (2019), 1305-1326.
 
 \vspace{1mm}
 
 \noindent
 [K] Carol R. Karp, 'Languages with expressions of infinite length',
 North-Holland Publishing Company, Amsterdam, (1964).
 
 \vspace{1mm}
 
 \noindent
 [Ko] Pascal Koiran, 'The theory of Liouville functions', The Journal of Symbolic Logic, 68(2), (2003), 353-365.
 \vspace{1mm}
 
 \noindent
 [MW] Angus Macintyre and A. J. Wilkie, 'On the decidability of the real exponential field', in 'Kreiseliana: about and around Georg Kreisel', edited by Piergiorgio Odifreddie, A. K. Peters (1996),441-467.

 \vspace{1mm}
 
 \noindent
 [PS1] Ya'acov Peterzil and Sergei Starchenko, 'Expansions of algebraically closed fields in o-minimal structures', Selecta Mathematica, New Series 7 (2001), 409-445.

 \vspace{1mm}
 
 \noindent
 [PS2] Ya'acov Peterzil and Sergei Starchenko, 'Expansions of algebraically closed fields II: functions of several variables', J. Math. Logic, Volume 3, Number 1 (2003), 1-35.
 
 \vspace{1mm}
 
 \noindent
 [PW] J. Pila and A. J. Wilkie, 'The rational points of a definable set', Duke Math. J. 133 (2006), no. 3, 591-616.
 
 \vspace{1mm}
 
 \noindent
 [R] David A. Ross, 'Yet another proof that the roots of a polynomial vary continuously as a function of the coefficients', arXiv:2207.00123v2. (2022).
 
 \vspace{1mm}
 
 \noindent
 [S] H. Sfouli, 'On the elementary theory of restricted real and imaginary parts of holomorphic functions', Notre Dame J. of Formal Logic, 53(1), (2012), 67-77.
 
 \vspace{1mm}
 
 \noindent
 [Sp] P. Speissegger,'Lectures on o-minimality', in 'Lectures on Algebraic Model Theory', edited by Bradd Hart and Matthew Valeriote, Fields Institute Monographs, AMS (2002), 47-65.
 
 \vspace{1mm}

 \noindent
 [T] M. Tsuji, 'Theory of meromorphic functions in a neighborhood of a closed set of capacity zero', Japanese J. Math. 19, (1944), 139-154.
 
 \vspace{1mm}
 
 \noindent
 [W1] A. J. Wilkie, 'Diophantine properties of sets definable in o-minimal structures', J. Symb. Logic, Vol 69, No. 3, (2004), 851-861.
 
 \vspace{1mm}

 \noindent
[W2] A. J. Wilkie, 'Liouville functions', in 'Logic Colloquium 2000: proceedings of the Annual European Summer Meeting of the Association of Symbolic Logic', edited by Ren\'{e} Cori, Alexander Razborov, Stevo Todor\u{c}evi\'{c} and Carol Wood, Lecture Notes in Logic 19, (2005), 383-391.

 \vspace{1mm}

 \noindent
 [W3] A. J. Wilkie, 'Some local definability theory for holomorphic functions', in 'Model Theory with Applications to Algebra and Analysis (Volume 1)', Volume 349 of LMS Lecture Note Series, CUP, (2008), 197-214. 

\vspace{1mm}

\noindent
[Z] Boris Zilber, 'Analytic and pseudo-analytic structures',  in 'Logic Colloquium 2000: proceedings of the Annual European Summer Meeting of the Association of Symbolic Logic', edited by Ren\'{e} Cori, Alexander Razborov, Stevo Todor\u{c}evi\'{c} and Carol Wood, Lecture Notes in Logic 19, (2005), 392-408.

 \end{document}